\documentclass{amsart}

\usepackage{microtype}
\usepackage{amssymb}
\usepackage{mathtools}
\usepackage{tikz-cd}
\usepackage{mathbbol} 

\usepackage[
	backend=biber,
	style=alphabetic,
	backref=true,
	url=false,
	doi=false,
	isbn=false,
	eprint=false]{biblatex}

\renewbibmacro{in:}{} 
\AtEveryBibitem{\clearfield{pages}} 

\DeclareFieldFormat{title}{\myhref{\mkbibemph{#1}}}
\DeclareFieldFormat
[article,inbook,incollection,inproceedings,patent,thesis,unpublished]
{title}{\myhref{\mkbibquote{#1\isdot}}}

\newcommand{\doiorurl}{%
	\iffieldundef{url}
	{\iffieldundef{eprint}
		{}
		{http://arxiv.org/abs/\strfield{eprint}}}
	{\strfield{url}}%
}

\newcommand{\myhref}[1]{%
	\ifboolexpr{%
		test {\ifhyperref}
		and
		not test {\iftoggle{bbx:eprint}}
		and
		not test {\iftoggle{bbx:url}}
	}
	{\href{\doiorurl}{#1}}
	{#1}%
}

\usepackage[
	bookmarks=true,
	linktocpage=true,
	bookmarksnumbered=true,
	breaklinks=true,
	pdfstartview=FitH,
	hyperfigures=false,
	plainpages=false,
	naturalnames=true,
	colorlinks=true,
	pagebackref=false,
	pdfpagelabels]{hyperref}

\hypersetup{
	colorlinks,
	citecolor=blue,
	filecolor=blue,
	linkcolor=blue,
	urlcolor=blue
}

\usepackage[capitalize, noabbrev]{cleveref}
\let\subsectionSymbol\S
\crefname{subsection}{\subsectionSymbol\!\!}{subsections}

\calclayout

\makeatletter
\@namedef{subjclassname@2020}{%
	\textup{2020} Mathematics Subject Classification}
\makeatother

\newtheorem{theorem}{Theorem}
\newtheorem*{theorem*}{Theorem}

\newtheorem*{proposition*}{Proposition}
\newtheorem{lemma}[theorem]{Lemma}
\newtheorem*{lemma*}{Lemma}

\newtheorem*{corollary*}{Corollary}
\theoremstyle{definition}

\newtheorem*{definition*}{Definition}
\newtheorem{remark}[theorem]{Remark}
\newtheorem*{remark*}{Remark}

\newtheorem*{example*}{Example}

\newtheorem*{construction*}{Construction}

\newtheorem*{convention*}{Convention}

\newtheorem*{terminology*}{Terminology}

\newtheorem*{notation*}{Notation}

\newtheorem*{question*}{Question}

\DeclareMathOperator{\bd}{\partial}
\DeclareMathOperator{\sign}{sign}
\newcommand{\ot}{\otimes}

\newcommand{\N}{\mathbb{N}}
\newcommand{\Z}{\mathbb{Z}}

\renewcommand{\k}{\Bbbk}
\newcommand{\sym}{\mathbb{S}}

\newcommand{\gsimplex}{\mathbb{\Delta}}
\newcommand{\gcube}{\mathbb{I}}

\newcommand{\Cat}{\mathsf{Cat}}
\newcommand{\Fun}{\mathsf{Fun}}
\newcommand{\Set}{\mathsf{Set}}
\newcommand{\Top}{\mathsf{Top}}

\newcommand{\Ch}{\mathsf{Ch}}
\newcommand{\simplex}{\triangle}
\newcommand{\sSet}{\mathsf{sSet}}
\newcommand{\cube}{\square}
\newcommand{\cSet}{\mathsf{cSet}}

\newcommand{\coAlg}{\mathsf{coAlg}}
\newcommand{\biAlg}{\mathsf{biAlg}}

\newcommand{\Mon}{\mathsf{Mon}}

\DeclareMathOperator{\forget}{U}
\DeclareMathOperator{\yoneda}{\mathcal{Y}}

\newcommand{\loops}{\Omega}
\DeclareMathOperator{\cobar}{\mathbf{\Omega}}

\DeclareMathOperator{\chains}{N}

\DeclarePairedDelimiter\bars{\lvert}{\rvert}

\DeclarePairedDelimiter\set{\{}{\}}

\newcommand{\id}{\mathsf{id}}
\renewcommand{\th}{\mathrm{th}}
\newcommand{\op}{\mathrm{op}}
\DeclareMathOperator*{\colim}{colim}

\newcommand{\Hom}{\mathrm{Hom}}

\newcommand{\coEnd}{\mathrm{coEnd}}

\newcommand{\xra}[1]{\xrightarrow{#1}}

\newcommand{\pdfEinfty}{\texorpdfstring{${E_\infty}$}{E-infty}}

\newcommand{\cC}{\mathcal{C}}

\newcommand{\cM}{\mathcal{M}}

\newcommand{\cO}{\mathcal{O}}
\newcommand{\cP}{\mathcal{P}}

\newcommand{\sC}{\mathsf{C}}

\newcommand{\fX}{\mathfrak{X}}

\newcommand{\bft}{\mathbf{t}}
\newcommand{\s}[1]{s^{\scalebox{.6}{{#1}}}}
\DeclareMathOperator{\aug}{\varepsilon}

\renewcommand{\S}{\mathbb{S}}
\newcommand{\biEnd}{\mathrm{biEnd}}

\newcommand{\As}{{\mathcal{A}\mathsf{s}}}
\newcommand{\Com}{\cC{om}}
\newcommand{\M}{\cM}
\newcommand{\UM}{{\forget(\M)}}

\DeclareMathOperator{\schains}{N^{\simplex}}
\DeclareMathOperator{\cchains}{N^{\cube}}

\DeclareMathOperator{\Schains}{S}
\DeclareMathOperator{\sSchains}{S^{\simplex}}
\DeclareMathOperator{\cSchains}{S^{\cube}}
\DeclareMathOperator{\schainsA}{N^{\simplex}_{\!\As}}
\DeclareMathOperator{\cchainsA}{N^{\cube}_{\!\As}}
\DeclareMathOperator{\SchainsA}{S_{\!\As}}
\DeclareMathOperator{\sSchainsA}{S^{\simplex}_{\!\As}}
\DeclareMathOperator{\cSchainsA}{S^{\cube}_{\!\As}}

\DeclareMathOperator{\cSchainsUM}{S^{\cube}_{\UM}}
\DeclareMathOperator{\schainsUM}{N^{\simplex}_{\UM}}
\DeclareMathOperator{\cchainsUM}{N^{\cube}_{\UM}}

\DeclareMathOperator{\sSing}{Sing^{\simplex}}
\DeclareMathOperator{\cSing}{Sing^{\cube}}

\DeclareMathOperator{\ccobar}{\mathbb{\Omega}}

\DeclareMathOperator{\ncobar}{\mathbb{\Omega}^{\mathrm{nec}}}

\newcommand{\Nec}{\mathsf{Nec}}
\newcommand{\nSet}{\mathsf{nSet}}

\newcommand{\pdfM}{\texorpdfstring{${\cM}$}{M}} 
\title[Adams' cobar construction as a monoidal $E_{\infty}$-coalgebra]{Adams' cobar construction as a monoidal $E_{\infty}$-coalgebra model of the based loop space}

\author[Medina-Mardones]{Anibal M. Medina-Mardones}
\address{Department of Mathematics, Western University, Canada}
\email{\href{anibal.medina.mardones@uwo.ca}{anibal.medina.mardones@uwo.ca}}

\author[Rivera]{Manuel Rivera}
\address{Department of Mathematics, Purdue University, USA}
\email{\href{mailto:manuelr@purdue.edu}{manuelr@purdue.edu}}


\usetikzlibrary{arrows,decorations.markings}
\tikzset{myptr/.style={decoration={markings,mark=at position 1 with %
			{\arrow[scale=2,>=stealth]{>}}},postaction={decorate}}}

\newsavebox\preproduct
\begin{lrbox}{\preproduct}
	\begin{tikzpicture}[scale=.3]
	\draw (0,0)--(0,-.8);
	\draw (0,0)--(.5,.5);
	\draw (0,0)--(-.5,.5);
	\end{tikzpicture}
\end{lrbox}
\newcommand{\product}{
	\usebox\preproduct}

\newsavebox\precoproduct
	\begin{lrbox}{\precoproduct}
		\begin{tikzpicture}[scale=.3]
		\draw (0,0)--(0,.8);
		\draw (0,0)--(.5,-.5);
		\draw (0,0)--(-.5,-.5);
		\end{tikzpicture}
	\end{lrbox}
\newcommand{\coproduct}{
	\usebox\precoproduct}

\newsavebox\preboundary
\begin{lrbox}{\preboundary}
	\begin{tikzpicture}[scale=.3]
	\draw (0,0)--(0,1.3);
	\draw (.5,0)--(.5,1.3);
	\draw [fill] (0,0) circle [radius=0.1];
	\draw (.9,.6)--(1.3,.6);
	\draw (1.7,0)--(1.7,1.3);
	\draw (2.2,0)--(2.2,1.3);
	\draw [fill] (2.2,0) circle [radius=0.1];
	\end{tikzpicture}
\end{lrbox}
\newcommand{\boundary}{
	\usebox\preboundary}

\newsavebox\preleftboundary
\begin{lrbox}{\preleftboundary}
	\begin{tikzpicture}[scale=.3]
	\draw (1.7,0)--(1.7,1.3);
	\draw (2.2,0)--(2.2,1.3);
	\draw [fill] (2.2,0) circle [radius=0.1];
	\end{tikzpicture}
\end{lrbox}
\newcommand{\leftboundary}{
	\usebox\preleftboundary}

\newsavebox\prerightboundary
\begin{lrbox}{\prerightboundary}
	\begin{tikzpicture}[scale=.3]
	\draw (0,0)--(0,1.3);
	\draw (.5,0)--(.5,1.3);
	\draw [fill] (0,0) circle [radius=0.1];
	\end{tikzpicture}
\end{lrbox}
\newcommand{\rightboundary}{
	\usebox\prerightboundary}

\newsavebox\precoboundary
\begin{lrbox}{\precoboundary}
	\begin{tikzpicture}[scale=.3]
	\draw (0,0)--(0,1.3);
	\draw (.5,0)--(.5,1.3);
	\draw [fill] (0,1.3) circle [radius=0.1];
	\draw (.9,.6)--(1.3,.6);
	\draw (1.7,0)--(1.7,1.3);
	\draw (2.2,0)--(2.2,1.3);
	\draw [fill] (2.2,1.3) circle [radius=0.1];
	\end{tikzpicture}
\end{lrbox}

\newsavebox\precounit
\begin{lrbox}{\precounit}
	\begin{tikzpicture}[scale=.3]
	\draw (0,0)--(0,1.3);
	\draw [fill] (0,0) circle [radius=0.1];
	\end{tikzpicture}
\end{lrbox}
\newcommand{\counit}{
	\usebox\precounit}

\newsavebox\preidentity
\begin{lrbox}{\preidentity}
	\begin{tikzpicture}[scale=.3]
	\draw (0,0)--(0,1.3);
	\end{tikzpicture}
\end{lrbox}

\newsavebox\preunit
\begin{lrbox}{\preunit}
	\begin{tikzpicture}[scale=.3]
	\draw (0,0)--(0,1.3);
	\draw [fill] (0,1.3) circle [radius=0.1];
	\end{tikzpicture}
\end{lrbox}

\newsavebox\preassociativity
\begin{lrbox}{\preassociativity}
	\begin{tikzpicture}[scale=.2]
	\path[draw] (0,0)--(0,1)--(-1,2);
	\draw (0,1)--(1,2);
	\draw (-.5,1.5)--(0,2);
	\draw (1.2,1)--(1.8,1);
	\path[draw] (3,0)--(3,1)--(2,2);
	\draw (3,1)--(4,2);
	\draw (3.5,1.5)--(3,2);
	\end{tikzpicture}
\end{lrbox}

\newsavebox\precoassociativity
\begin{lrbox}{\precoassociativity}
	\begin{tikzpicture}[scale=.2]
	\path[draw] (0,0)--(0,-1)--(-1,-2);
	\draw (0,-1)--(1,-2);
	\draw (-.5,-1.5)--(0,-2);
	\draw (1.2,-1)--(1.8,-1);
	\path[draw] (3,0)--(3,-1)--(2,-2);
	\draw (3,-1)--(4,-2);
	\draw (3.5,-1.5)--(3,-2);
	\end{tikzpicture}
\end{lrbox}
\newcommand{\coassociativity}{
	\usebox\precoassociativity}

\newsavebox\preleftcomb
\begin{lrbox}{\preleftcomb}
	\begin{tikzpicture}[scale=.2]
	\path[draw] (0,0)--(0,-1)--(-1,-2);
	\draw (0,-1)--(1,-2);
	\draw (-.5,-1.5)--(0,-2);
	\end{tikzpicture}
\end{lrbox}

\newsavebox\prerightcomb
\begin{lrbox}{\prerightcomb}
	\begin{tikzpicture}[scale=.2]
	\path[draw] (3,0)--(3,-1)--(2,-2);
	\draw (3,-1)--(4,-2);
	\draw (3.5,-1.5)--(3,-2);
	\end{tikzpicture}
\end{lrbox}

\newsavebox\preinvolution
\begin{lrbox}{\preinvolution}
	\begin{tikzpicture}[scale=.2]
	\path[draw] (0,0)--(0,.5)--(-.5,1)--(0,1.5)--(0,2);
	\path[draw] (0,.5)--(.5,1)--(0,1.5);
	\end{tikzpicture}
\end{lrbox}

\newsavebox\preleftcounitality
\begin{lrbox}{\preleftcounitality}
	\begin{tikzpicture}[scale=.3]
	\draw (0,0)--(0,.8);
	\draw (0,0)--(.5,-.5);
	\draw (0,0)--(-.5,-.5);
	\draw [fill] (-.5,-.5) circle [radius=0.1];
	\draw (.7,0)--(1.1,0);
	\path[draw] (1.5,-.5)--(1.5,.8);
	\end{tikzpicture}
\end{lrbox}
\newcommand{\leftcounitality}{
	\usebox\preleftcounitality}

\newsavebox\preleftcounitcoproduct
\begin{lrbox}{\preleftcounitcoproduct}
	\begin{tikzpicture}[scale=.3]
	\draw (0,0)--(0,.8);
	\draw (0,0)--(.5,-.5);
	\draw (0,0)--(-.5,-.5);
	\draw [fill] (-.5,-.5) circle [radius=0.1];
	\end{tikzpicture}
\end{lrbox}

\newsavebox\prerightcounitality
\begin{lrbox}{\prerightcounitality}
	\begin{tikzpicture}[scale=.3]
	\draw (0,0)--(0,.8);
	\draw (0,0)--(.5,-.5);
	\draw (0,0)--(-.5,-.5);
	\draw [fill] (.5,-.5) circle [radius=0.1];
	\draw (-.7,0)--(-1.1,0);
	\path[draw] (-1.5,-.5)--(-1.5,.8);
	\end{tikzpicture}
\end{lrbox}
\newcommand{\rightcounitality}{
	\usebox\prerightcounitality}

\newsavebox\prerightcounitcoproduct
\begin{lrbox}{\prerightcounitcoproduct}
	\begin{tikzpicture}[scale=.3]
	\draw (0,0)--(0,.8);
	\draw (0,0)--(.5,-.5);
	\draw (0,0)--(-.5,-.5);
	\draw [fill] (.5,-.5) circle [radius=0.1];
	\end{tikzpicture}
\end{lrbox}

\newsavebox\preleftunitality
\begin{lrbox}{\preleftunitality}
	\begin{tikzpicture}[scale=.3]
	\draw (0,0)--(0,-.8);
	\draw (0,0)--(-.5,.5);
	\draw (0,0)--(.5,.5);
	\draw [fill] (.5,.5) circle [radius=0.1];
	\draw (.7,0)--(1.1,0);
	\path[draw] (1.5,.5)--(1.5,-.8);
	\end{tikzpicture}
\end{lrbox}

\newsavebox\prerightunitality
\begin{lrbox}{\prerightunitality}
	\begin{tikzpicture}[scale=.3]
	\draw (0,0)--(0,-.8);
	\draw (0,0)--(-.5,.5);
	\draw (0,0)--(.5,.5);
	\draw [fill] (-.5,.5) circle [radius=0.1];
	\draw (-.7,0)--(-1.1,0);
	\path[draw] (-1.5,.5)--(-1.5,-.8);
	\end{tikzpicture}
\end{lrbox}

\newsavebox\preproductcounit
\begin{lrbox}{\preproductcounit}
	\begin{tikzpicture}[scale=.3]
	\draw (0,0)--(0,-.8);
	\draw (0,0)--(.5,.5);
	\draw (0,0)--(-.5,.5);
	\draw [fill] (0,-.8) circle [radius=0.1];
	\end{tikzpicture}
\end{lrbox}
\newcommand{\productcounit}{
	\usebox\preproductcounit}

\newsavebox\preunitcoproduct
\begin{lrbox}{\preunitcoproduct}
	\begin{tikzpicture}[scale=.3]
	\draw (0,0)--(0,.8);
	\draw (0,0)--(.5,-.5);
	\draw (0,0)--(-.5,-.5);
	\draw [fill] (0,.8) circle [radius=0.1];
	\end{tikzpicture}
\end{lrbox}

\newsavebox\preleibniz
\begin{lrbox}{\preleibniz}
	\begin{tikzpicture}[scale=.245]
	\draw (0,.3)--(0,-.3);
	\draw (0,.3)--(.5,.8);
	\draw (0,.3)--(-.5,.8);
	\draw (0,-.3)--(0,.3);
	\draw (0,-.3)--(.5,-.8);
	\draw (0,-.3)--(-.5,-.8);

	\draw (.7,0)--(1.1,0);
	\draw (2.1,.8)--(2.1,.3)--(1.7,0)--(1.7,-.8);
	\draw (2.1,.3)--(2.9,-.3);
	\draw (3.3,.8)--(3.3,0)--(2.9,-.3)--(2.9,-.8);

	\draw (3.8,0)--(4.1,0);
	\draw (4.6,.8)--(4.6,0)--(5,-.3)--(5,-.8);
	\draw (5,-.3)--(5.8,.3);
	\draw (5.8,.8)--(5.8,.3)--(6.2,0)--(6.2,-.8);
	\end{tikzpicture}
\end{lrbox}

\newsavebox\prebialgebra
\begin{lrbox}{\prebialgebra}
	\begin{tikzpicture}[scale=.5]
	\draw (0,.3)--(0,-.3);
	\draw (0,.3)--(.5,.8);
	\draw (0,.3)--(-.5,.8);
	\draw (0,-.3)--(0,.3);
	\draw (0,-.3)--(.5,-.8);
	\draw (0,-.3)--(-.5,-.8);

	\draw (.8,0)--(1.2,0);

	\draw (2.1,.8)--(2.1,.3)--(1.7,0)--(2.1,-.3)--(2.1,-.8);

	\draw (3.1,.8)--(3.1,.3)--(3.5,0)--(3.1,-.3)--(3.1,-.8);

	\draw (2.1,.3)--(3.1,-.3);
	\draw (2.1,-.3)--(2.5,-.06);
	\draw (3.1,.3)--(2.7,.06);
	\end{tikzpicture}
\end{lrbox}

\newsavebox\precommutativity
\begin{lrbox}{\precommutativity}
	\begin{tikzpicture}[scale=.27]
	\draw (.3,0)--(.3,-1);
	\draw (.3,0)--(.8,.5);
	\draw (.3,0)--(-.2,.5);

	\draw (1.2,-.4)--(1.8,-.4);

	\draw (3,-.3)--(3,-1);
	\draw (2.5,0)--(3,-.3);
	\draw (3.5,0)--(3,-.3);
	\draw (2.46,0)--(2.9,.18);
	\draw (3.1,.3)--(3.5,.5);
	\draw (3.5,0)--(2.5,.5);
		\end{tikzpicture}
	\end{lrbox}

\newsavebox\precocommutativity
\begin{lrbox}{\precocommutativity}
	\begin{tikzpicture}[scale=.3]
		\draw (0,0)--(0,.8);
		\draw (0,0)--(.5,-.5);
		\draw (0,0)--(-.5,-.5);
		\draw (.7,0)--(1.1,0);
		\node[scale=.8] at (1.6,0){$\sigma$};
		\draw (2.2,0)--(2.2,.8);
		\draw (2.2,0)--(2.7,-.5);
		\draw (2.2,0)--(1.7,-.5);
	\end{tikzpicture}
\end{lrbox}
\newcommand{\cocommutativity}{
	\usebox\precocommutativity}

\newsavebox\cellprecoboundary
\begin{lrbox}{\cellprecoboundary}
	\begin{tikzpicture}[scale=.3]
		\draw (-2,0)--(-2,1.3);
		\draw (-2.5,0)--(-2.5,1.3);
		\draw [fill] (-2.5,0) circle [radius=0.1];
		\draw (-2,-.5)--(2,-.5); 
		\draw (2,0)--(2,1.3);
		\draw (2.5,0)--(2.5,1.3);
		\draw [fill] (2.5,0) circle [radius=0.1];
		\draw (0,.7)--(0,0);
		\draw (0,.7)--(.5,1.3);
		\draw (0,.7)--(-.5,1.3);
		\node at (0,-.6){};
	\end{tikzpicture}
\end{lrbox}
\newcommand{\cellcoboundary}{
	\usebox\cellprecoboundary}

\newsavebox\cellprecoproduct
\begin{lrbox}{\cellprecoproduct}
	\begin{tikzpicture}[scale=.3]
		\draw (0,0)--(0,.8);
		\draw (0,0)--(.5,-.5);
		\draw (0,0)--(-.5,-.5);
		\draw [fill] (0,-.9) circle [radius=0.05];
		\node at (0,-1){};
	\end{tikzpicture}
\end{lrbox}
\newcommand{\cellcoproduct}{
	\usebox\cellprecoproduct}

\newsavebox\cellprecounit
\begin{lrbox}{\cellprecounit}
	\begin{tikzpicture}[scale=.3]
		\draw (0,0)--(0,1.3);
		\draw [fill] (0,0) circle [radius=0.1];
		\draw [fill] (0,-.5) circle [radius=0.05];
		\node at (0,-.6){};
	\end{tikzpicture}
\end{lrbox}
\newcommand{\cellcounit}{
	\usebox\cellprecounit}

\begin{document}
	\vspace*{-1cm}

\begin{abstract}
	We prove that the classical map comparing Adams' cobar construction on the singular chains of a pointed space and the singular cubical chains on its based loop space is a quasi-isomorphism preserving explicitly defined monoidal $E_\infty$-coalgebra structures. This contribution extends to its ultimate conclusion a result of Baues, stating that Adams' map preserves monoidal coalgebra structures.
\end{abstract}
	\maketitle
	\vspace*{-1.2cm}
	\tableofcontents

\section{Introduction}

For any topological space $\fX$, its complex of simplicial or cubical singular chains $\Schains(\fX)$ -- regarded as a differential graded (dg) abelian group -- encodes the homology of $\fX$ in its quasi-isomorphism type.
More homotopical information can be stored in the quasi-isomorphism type of this chain complex if considered as a (coassociative) coalgebra, which we will denote $\SchainsA(\fX)$, where the coproduct comes from a natural choice of chain approximation to the diagonal $\fX \to \fX \times \fX$.
For instance, the cohomology ring of $\fX$ is retained, but the action of the Steenrod algebra on its mod~$p$ cohomology is not.

In Mandell's seminal work \cite{mandell2006homotopy_type} it is shown that, when $\fX$ is nilpotent and finite type, the entire homotopy type of $\fX$ can be encoded in the quasi-isomorphism type of this complex if considered as an $E_\infty$-coalgebra, a structure providing $\SchainsA(\fX)$ with coherent homotopies witnessing the derived cocommutativity of the coproduct coming from the strict symmetry of the diagonal map.

The first contribution of this paper is to explicitly endow the cubical singular chains of the based loop space $\loops_x \fX$, with the structure of a monoidal $E_\infty$-coalgebra extending the Serre diagonal.
More specifically, we verify that the monoid structure induced on $\cSchains(\loops_x \fX)$ by the concatenation of loops is compatible with a natural $E_\infty$-coalgebra structure on cubical singular chains, similar to the one defined in \cite{medina2022cube_einfty}.

Applying Adams' cobar construction to the coalgebra of simplicial singular chains of $(\fX, x)$, one obtains another monoidal algebraic model $\cobar \sSchainsA(\fX, x)$ of $\loops_x \fX$ \cite{adams1956cobar}.
More precisely, Adams constructed a natural monoidal chain map $\theta$ from $\cobar \sSchainsA(\fX, x)$ to $\cSchains(\loops_x \fX)$ and proved it to be a quasi-isomorphism if $\fX$ is simply-connected, a statement that also holds true for path-connected spaces after \cite{rivera2018cubical}.
The model $\cobar \sSchainsA(\fX, x)$ is smaller than $\cSchains(\loops_x \fX)$ and unlocks effective analysis of quantitative and qualitative properties of $\loops_x \fX$, as illustrated for instance in \cite{chainalgebraloops} and \cite{adamscobarequivalence}.

The second main contribution of this paper is to make Adams model into a monoidal $E_\infty$-coalgebra and to prove that
\[
\theta \colon \cobar \sSchainsA(\fX, x) \to \cSchains(\loops_x \fX)
\]
respects this higher structure.
Although not pursued in the present article, we remark that the explicit nature of our $E_\infty$-extension invites the study of primary and secondary operations for loops spaces using Adams' model and the tools developed in \cite{medina2021may_st}, \cite{medina2020cartan}, \cite{medina2021adem}, and \cite{medina2021comch}.

Our starting point is groundbreaking work by Baues, which imply statements similar to those in this work but in the category of (coassociative) coalgebras.
Baues reinterpreted Adams' algebraic construction at a deeper geometric level \cite{baues1998hopf}, which allowed him to endow $\cobar \sSchainsA(\fX, x)$ with the structure of a monoidal coalgebra, and to show that $\theta$ preserves this structure.
To prove our statement we interpret Adams' construction at an even deeper categorical level.
We interpret Baues' geometric cobar construction, originally defined for $1$-reduced simplicial sets, as a functor
\begin{equation*}
	\ccobar \colon \sSet^0 \to \Mon_{\cSet},
\end{equation*}
from the category of $0$-reduced simplicial sets to that of monoidal cubical sets.
The key difference with Baues' original work is the use of connections to obtain a natural construction before geometric realization.

Additionally, we need a suitable model of the $E_\infty$-operad endowing cubical chains with a natural $E_\infty$-coalgebra extending the Serre diagonal.
For this we take the operad $\UM$ introduced in \cite{medina2020prop1}.
After proving that its coalgebras form a monoidal category, we show that the functor $\cchainsUM \colon \cSet \to \coAlg_\UM$ -- defined in \cite{medina2022cube_einfty} with a different sign convention -- is monoidal.
This allows us to construct the following extension of Adams and Baues' structures.

\begin{theorem*}
	The following diagram commutes up to natural isomorphisms:
	\[
	\begin{tikzcd} [row sep=small]
		& \Mon_{\coAlg_\UM} \arrow[d] \\
		\Mon_{\cSet} \arrow[ru, "\cchainsUM", out=70, in=180, near start] \arrow[r, "\cchainsA"]
		& \Mon_{\coAlg} \arrow[d] \\
		\sSet^0 \arrow[r, "\cobar \schainsA"] \arrow[u, "\ccobar"]
		& \Mon_{\Ch},
	\end{tikzcd}
	\]
	where the unlabeled arrows are forgetful functors.
\end{theorem*}

In the diagram of the above theorem, the arrow from $\sSet^0$ to $\Mon_{\Ch}$ is Adams' cobar construction, the one from $\sSet^0$ to $\Mon_{\coAlg}$ is Baues' enhancement, and the one from $\sSet^0$ to $\Mon_{\coAlg_{\UM}}$ is our lift.
Additionally, we prove the following statement about Adams's map.

\begin{theorem*}
	For any pointed space $(\fX, x)$,
	\[
	\theta \colon \cobar \sSchainsA(\fX, x) \to \cSchains(\loops_x \fX)
	\]
	is a quasi-isomorphism of monoidal $\UM$-coalgebras.
\end{theorem*}

The fact that $\theta$ respects the monoid structure in $\Ch$ was proven by Adams, whereas the compatibility of the monoid structure with the Serre coalgebra structure was established by Baues.
Our contribution is the compatibility of the monoid structure with a full $E_\infty$-coalgebra extension of Serre's coalgebra.
We also remark that, whereas both Adams and Baues worked in the setting where the underlying space is simply connected, the above theorem does not require any connectivity or finiteness hypotheses.

\subsection*{Related work}

Kadeishvili \cite{kadeishvili1999coproducts, kadeishvili2003cupi} explicitly described monoidal cup-$i$ coproducts on $\cobar \schainsA(X)$ extending Baues coalgebra.
Kadeishvili, as Baues, used cubical methods to define these coproducts and to compare them, in the $1$-connected setting, to cup-$i$ coproducts extending the Serre coalgebra structure on the cubical singular chains of the based loop space.
Additionally, there are several papers \cite{smirnov1990iterated, smith1994cobar, smith2000operads, kadeishvili1998iterating} that predict the existence of, but do not construct, an $E_\infty$-structure on the cobar construction on the chains of simply connected simplicial sets.

On the dual side, Fresse \cite{fresse2003hopf} provided the bar construction of an algebra over the surjection operad with the structure of a comonoid in the category of algebras over the Barratt--Eccles operad.
Additionally, in \cite{fresse2010bar} he used a model category structure on reduced operads \cite{berger2003modelcategory, hinich1997homologicalalgebra} to iterate the bar construction on algebras over cofibrant $E_\infty$-operads.

The use of coalgebras instead of algebras allows us to relate the cobar construction to the based loop space directly --via the Adams map-- without imposing restrictions on the underlying homotopy type, as done by Fresse.
Furthermore, by grounding our approach on the cubical perspective at the heart of Adams' and Baues' seminal papers, we are able to preserve the natural monoidal structures when defining our $E_\infty$-enhancements.
	
\section*{Acknowledgements}

We would like to thank Clemens Berger, Matthias Franz, Kathryn Hess, Ralph Kaufmann, Emilio Minichello, Viet-Cuong Pham, Dennis Sullivan, and Mahmoud Zeinalian, for insightful discussions related to this project.

Both authors are grateful for the support and excellent working conditions of the Max Planck Institute for Mathematics in Bonn, Germany.
A.M. acknowledges support from \textit{Innosuisse grant} \mbox{32875.1 IP-ICT-1} and ANR 20 CE40 0016 01 PROJET HighAGT.
M.R. acknowledges support from \textit{NSF Grant DMS 210554} and the \textit{Karen EDGE Fellowship}.

\section{Conventions and preliminaries}\label{s:preliminaries}

\subsection{Coalgebras}\label{ss:coalgebras}

Throughout this article $\k$ denotes a commutative and unital ring and we work over its associated symmetric monoidal category of (homologically) graded chain complexes of $\k$-modules $(\Ch, \ot, \k)$.

A \textit{coalgebra} consists of a chain complex $C$ and chain maps $\Delta \colon C \to C \ot C$ and $\varepsilon \colon C \to \k$ satisfying the usual coassociativity and counitality relations.
Denote by $\coAlg$ the category of coalgebras with morphisms being structure preserving chain maps.
The category $\coAlg$ is symmetric monoidal, with braiding induced from $\Ch$ and structure maps of a product $C \ot C^\prime$ given by
\begin{gather*}
	C \ot C^\prime \xra{\Delta \ot \Delta^\prime}
	(C \ot C) \ot (C^\prime \ot C^\prime) \xra{(23)}
	(C \ot C^\prime) \ot (C \ot C^\prime), \\
	C \ot C^\prime \xra{\varepsilon \ot \varepsilon^\prime}
	\k \ot \k \xra{\cong} \k.
\end{gather*}

A \textit{coaugmentation} on a coalgebra $C$ is a coalgebra map $\nu \colon \k \to C$.
A coalgebra is said to be \textit{coaugmented} if it is equipped with a coaugmentation.
We denote by $\coAlg^\ast$ the category of coaugmented coalgebras with morphisms being coaugmentation preserving coalgebra maps.
A coaugmented coalgebra is a \textit{connected coalgebra} if it is $0$ in negative degrees and the coaugmentation induces an isomorphism $\k \cong C_0$ of $\k$-modules.
We denote by $\coAlg^0$ the full subcategory of $\coAlg^\ast$ defined by these.

\subsection{Monoids}

A \textit{monoidal object} in a monoidal category $(\sC, \ot, \mathbb{1})$ is an object $M$ together with morphisms $\mu \colon M \ot M \to M$ and $\eta \colon \mathbb{1} \to M$ satisfying the usual associativity and unital relations.
The category of these together with structure preserving morphisms, referred to as \textit{monoidal morphisms}, is denoted $\Mon_\sC$.
We remark that a lax monoidal functor $\sC \to \sC^\prime$ induces a functor between their categories of monoids $\Mon_\sC \to \Mon_{\sC^\prime}$.
For example, the forgetful functor from $\coAlg$ to $\Ch$ is monoidal and so, it induces a forgetful functor from monoidal coalgebras to monoidal chain complexes, which are more commonly known as bialgebras and algebras respectively, but this terminology is not well suited for our purposes.

\subsection{Simplicial theory}\label{ss:simplicial}

The \textit{simplex category} is denoted by $\simplex$ and its objects by $[n]$.
The morphisms in $\simplex$ are generated by coface maps, denoted by $\partial^{j}_n \colon [n-1] \to [n]$ for $0\leq j \leq n$, and codegeneracy maps, denoted by $\xi^{j}_n \colon [n+1] \to [n]$ for $0 \leq j \leq n+1$.
These satisfy the usual cosimplicial identities.
For simplicity, we omit the subscript in the notation when there is no risk of confusion and simply denote these maps by $\partial^{j} \colon [n-1] \to [n]$ and $\xi^{j} \colon [n+1] \to [n]$.

The category of \textit{simplicial sets} $\Fun(\simplex^\op, \Set)$ is denoted by $\sSet$ and the \textit{standard $n$-simplex} $\simplex(-, [n])$ by $\simplex^n$.
For any simplicial set $X$ we write, as usual, $X_n$ instead of $X[n]$, and identify the elements of $\simplex^n_m$ with increasing tuples $[v_0, \dots, v_m]$ where $v_i \in \{0, \dots, n\}$.

If $X$ is such that $X_0$ is a singleton set we say that $X$ is \textit{reduced}.
We denote the full subcategory of reduced simplicial sets by $\sSet^0$.

We consider the \textit{topological $n$-simplex} $\gsimplex^n$ embedded in $\mathbb{R}^{n+1}$ as
\[ \gsimplex^n = \{ (x_0, \dots, x_n) \in \mathbb{R}^{n+1} | x_0^2 + \dots + x_n^2=1\} \]
so that the $i$-th vertex is given by the standard basis vector $v_i=(0,\ldots,1, \ldots,0) \in \mathbb{R}^{n+1}$ with $1$ in the $i$-th entry and $0$ in all other entries.
Consider the usual adjunction pair formed by the \textit{geometric realization} and \textit{singular complex}
\[
\begin{tikzcd}
	\bars{\ } \colon \sSet \arrow[r, shift left=.5ex] &
	\Top :\! \sSing, \arrow[l, shift left=.5ex]
\end{tikzcd}
\]
determined by the spaces $\gsimplex^n$ and usual coface inclusions $\delta^i \colon \gsimplex^{n-1} \to \gsimplex^n$ and codegeneracy projections $\sigma^j \colon \gsimplex^n \to \gsimplex^{n-1}$.

The functor of (normalized) \textit{simplicial chains} $\schains \colon \sSet \to \Ch$ is defined on any simplicial set $X$ by first considering the chain complex $(\k[X], \partial)$, given on degree $n$ by the free $\k$-module generated by the $n$-simplices of $X$ with differential given by the alternating sum of face maps, and then modding out by the sub-chain complex of degenerate elements.
We remark that this functor is naturally equivalent to the composition of the geometric realization functor and that of cellular chains with respect to the standard CW structure.
When no confusion arises from doing so we write $\chains$ instead of $\schains$ and refer to it simply as the functor of \textit{chains}.
We will denote the functor of (simplicial) \textit{singular chains} $\schains \circ \sSing \colon \Top \to \Ch$ by $\sSchains$, where $\sSing \colon \Top \to \sSet$ denotes the singular complex functor. We modify this construction for a pointed topological space $(\fX,x)$ by only considering maps $\gsimplex^n \to \fX$ sending all vertices to $x$.
This produces a reduced simplicial set $\sSing(\fX,x)$ whose normalized chains we denote by $\sSchains(\fX,x)$.

We now recall a classical chain approximation of the diagonal map on the chains of a simplicial set, i.e. a lift of the functor of chains to coalgebras:
\[
\begin{tikzcd}
	& \coAlg \arrow[d] \\
	\sSet \arrow[r, "\schains"] \arrow[ur, "\schainsA", out=70, in=180] & \Ch.
\end{tikzcd}
\]
It suffices to define this functor $\schainsA$ on standard simplices. For any $n \in \N$, define $\Delta \colon \chains(\simplex^n) \to \chains(\simplex^n)^{\ot2}$ by
\[
\Delta \big( [v_0, \dots, v_q] \big) = \sum_{i=0}^q \ [v_0, \dots, v_i] \ot [v_i, \dots, v_q]
\]
and $\epsilon \colon \chains(\simplex^n) \to \k$ by
\[
\epsilon \big( [v_0, \dots, v_q] \big) = \begin{cases} 1 & \text{ if } q = 0, \\ 0 & \text{ if } q>0. \end{cases}
\]
We referred to this lift as the \textit{Alexander--Whitney coalgebra} structure on simplicial chains.

If $X$ is a reduced simplicial set then $\schainsA(X)$ becomes a connected coalgebra with the coaugmentation $\nu \colon \k \to \schains(X)$ induced by sending $1$ to the basis element represented by the unique $0$-simplex of $X$.
Hence, the functor $\schainsA$ restricts to a functor
\[
\schainsA \colon \sSet^0 \to \coAlg^0.
\]

\subsection{Cubical theory}\label{ss:cubical}

We include an abridged presentation of cubical sets parallel to the our treatment of simplicial sets in the previous section.
For more details we refer the reader to any of \cite{grandis2003cubical, medina2022cube_einfty, medina2023flowing}.

The \textit{cube category} (with a connection) $\cube$ is the subcategory of $\Set$ whose objects are $2^n = \{0, 1\}^n$ for $n \in \N$ with $2^0 = \{0\}$.
We will denote $2^1$ by $2$ and $2^0$ by $1$, which is the unit of the monoidal structure of $\cube$ given by $2^n \times 2^m = 2^{n+m}$.
The morphisms of $\cube$ are monoidally generated by
\begin{equation}\label{eq:cube generators}
	\begin{tikzcd}
		1 \arrow[r, out=45, in=135, "\delta_0"] \arrow[r, out=-45, in=-135, "\delta_1"'] & 2 \arrow[l, "\sigma"'] &[-10pt] \arrow[l, "\gamma"'] 2 \times 2
	\end{tikzcd}
\end{equation}
defined by
\begin{gather*}
	\delta_0(0) = 0, \qquad
	\delta_1(0) = 1, \qquad
	\sigma(p) = 0, \qquad
	\gamma(p,q) = \max(p,q),
\end{gather*}
for $p,q \in \{0,1\}$.
We remark that the cube category without connections will not be considered in this work.

The category of \textit{cubical sets} $\Fun(\cube^\op, \Set)$ is denoted by $\cSet$ and
the \textit{standard $n$-cube} $\cube(-, 2^n)$ by $\cube^n$.
For any cubical set $Y$ we write, as usual, $Y_n$ instead of $Y(2^n)$.
The monoidal structure on $\cube$ induces one on $\cSet$.
More explicitly, for two cubical sets we have
\[
(Y \times Y')_n = \bigsqcup_{i+j=n} Y_i \times Y'_j.
\]

Consider the monoidal functor $\cube \to \Top$ defined by assigning $2$ to the usual interval $\gcube^1$ and \eqref{eq:cube generators} to the continuous maps
\[
\begin{tikzcd}
	\gcube^0 \arrow[r, out=45, in=135, "\delta_0"] \arrow[r, out=-45, in=-135, "\delta_1"'] & \gcube^1 \arrow[l, "\sigma"'] &[-10pt] \arrow[l, "\gamma"'] \gcube^2,
\end{tikzcd}
\]
where $\delta_i$ and $\sigma$ are the canonical cubical co-face and co-degeneracy maps, respectively, and $\gamma(s,t)=\text{max}(s,t)$.

From this functor we obtain the usual adjunction pair formed by the \textit{geometric realization} and \textit{singular complex}
\[
\begin{tikzcd}
	\bars{\ } \colon \cSet \arrow[r, shift left=.5ex] &
	\Top :\! \cSing. \arrow[l, shift left=.5ex]
\end{tikzcd}
\]
Notice that $\cSing$ is lax monoidal
\[
(\gcube^n \to \fX) \times (\gcube^m \to \fX') \mapsto (\gcube^{n+m} \to \fX \times \fX').
\]

The functor of (normalized) \textit{cubical chains} $\cchains \colon \cSet \to \Ch$ is the unique monoidal functor defined by assigning to $\cube^1$ the usual cellular chains of the interval $\gcube^1$.
Explicitly, it is defined on any cubical set $Y$ by first considering the chain complex $(\k[Y], \delta)$, given on degree $n$ by the free $\k$-module generated by the $n$-cubes of $Y$ with differential $\delta$ given on any $n$-cube $y \in Y_n$ by
\[
\delta(y) = \sum_{i=1}^n (-1)^i
\big(
Y(\text{id}_2^{i-1} \times \delta_1 \times \text{id}_2^{n-i})(y) -
Y(\text{id}_2^{i-1} \times \delta_0 \times \text{id}_2^{n-i})(y)
\big),
\]
and then modding out by the sub-complex of degenerate cubes.
When no confusion arises from doing so we write $\chains$ instead of $\cchains$ and refer to it simply as the functor of \textit{chains}.
We remark that $\cchains$ is naturally equivalent to the composition of the geometric realization and cellular chains functors with respect to the canonical CW structure.
We will denote the functor of (cubical) \textit{singular chains} $\cchains \circ \cSing \colon \Top \to \Ch$ by $\cSchains$.

Since $\chains(\cube^1)$ is a coalgebra and the category of coalgebras is monoidal, there is a unique monoidal functor $\cchainsA$ lifting the functor of chains
\[
\begin{tikzcd}
	& \coAlg \arrow[d] \\
	\cSet \arrow[r, "\cchains"] \arrow[ur, "\cchainsA", out=70, in=180] & \Ch.
\end{tikzcd}
\]
We refer to this lift as the \textit{Serre coalgebra} structure on cubical chains, and to the coproduct $\Delta$ as the Serre diagonal.

\section{Adams's model of the based loop space}\label{s:theorem1}

In this section we revisit Adams's classical cobar construction as a model for the based loop space.
A deeper exploration of Adams's comparison map from the cobar construction of the simplicial singular chains on a space to the cubical chains on the based loop space naturally leads us to the framework of necklical sets, a notion related to both simplicial sets and cubical sets.
We use this framework to construct a functorial cubical model for the based loop space.
Similar constructions and results may be found in \cite{baues1980geometry}, \cite{berger1995loops}, \cite{baues1998hopf}, \cite{dugger2011rigidification}, \cite{galvez2020hopf}, and \cite{rivera2018cubical, rivera2019path}.

\subsection{The based loop space}

The \textit{path space} functor
\[
P \colon \Top \to \Top
\]
assigns to $\fX$ the space
\[
P(\fX) = \big\{\alpha \colon [0,r] \to \fX \mid \text{$\alpha$ continuous and } r \in [0,\infty)\big\}
\]
equipped with the compact-open topology.
For any $x,x' \in \fX$ the subspace $P(\fX;x,x')$ consists of those paths starting and ending at $x$ and $x'$ respectively. The points of $P(\fX)$ are often referred to as \textit{Moore paths on $\fX$}.
We will think of this construction as a functor on bipointed spaces in the obvious way.
There is a composition structure
\begin{equation}\label{eq:composition structure}
	\begin{split}
		P(\fX;x,x') \times P(\fX;x',x'') &\to P(\fX;x,x'') \\
		(\alpha, \beta) &\mapsto \alpha \cdot \beta
	\end{split}
\end{equation}
given by concatenation of paths with addition of parameters, whose identities are constant paths with $r=0$.
We may assemble this data into a topologically enriched category $\mathcal{P}(\fX)$ that has the points of $\fX$ as objects, the spaces $P(\fX;x,x')$ as morphisms from $x$ to $x'$, concatenation of paths as composition, and constant paths with $r=0$ as identity morphisms.

Denote by $\Top^\ast$ the category of pointed topological spaces, which we think of as a full subcategory of that of bipointed spaces in the obvious way.
The \textit{based loop space} functor
\[
\loops \colon \Top^\ast \to \Mon_{\Top}
\]
associates to a pointed space $(\fX,x)$ the space
\[
\loops_x \fX = P(\fX;x,x)
\]
of loops in $\fX$ based at $x$ with the monoid structure induced from the composition structure \eqref{eq:composition structure}.

Since $\cSing$ is lax monoidal, $\cSing(\loops_x \fX)$ is a monoid in $\cSet$, and, since $\cchains$ is monoidal, $\cSchains(\fX,x)$ is a monoid in $\Ch$.

\subsection{The cobar construction}\label{ss:cobar construction}

We now describe an algebraic analogue of the based loop space introduced by Adams \cite{adams1956cobar}.
The \textit{cobar construction} is the functor
\[
\cobar \colon \coAlg^\ast \to \Mon_{\Ch}
\]
defined on objects as follows.
Let $(C, \Delta, \varepsilon, \nu)$ be a coaugmented coalgebra.
Denote by $\overline{C}$ the cokernel of the coaugmentation $\nu \colon \k \to C$ and recall that $\s{}$ is the suspension functor.
The cobar construction $\cobar C$ of this coaugmented coalgebra is the graded module
\[
T(\s{-1} \overline{C}) =
\k \oplus \s{-1}\overline{C} \oplus (\s{-1}\overline{C})^{\ot 2} \oplus (\s{-1}\overline{C})^{\ot 3} \oplus \dots
\]
regarded as a monoid in $\Ch$ with free associative product $\mu \colon T(\s{-1} \overline{C})^{\ot 2} \to T(\s{-1} \overline{C})$ given by concatenation, unit map $\eta \colon \k \to T(\s{-1} \overline{C})$ by the obvious inclusion, and differential constructed by extending the linear map
\[
-\, \s{-1} \circ \partial \circ \s{+1} \ + \ (\s{-1} \ot \s{-1}) \circ \Delta \circ \s{+1} \ \colon \
\s{-1} \overline{C} \to \s{-1}\overline{C} \oplus (\s{-1}\overline{C} \ot \s{-1}\overline{C}) \hookrightarrow T(\s{-1}C)
\]
as a derivation using the freeness of the underlying graded monoid.
On morphisms, the functor $\cobar$ is defined using the functoriality of the free graded monoid construction.
For any $x_1, \dots, x_k \in \overline{C}$, we denote
\[
[x_1| \cdots | x_k]= \s{-1} x_1 \otimes \cdots \otimes \s{-1}x_k \in (\s{-1}\overline{C})^{\otimes k}.
\]

\subsection{Adams's map}\label{ss:adams maps}

Adams made precise the sense in which the cobar functor may be understood as an algebraic analogue of the based loop space functor.
This was achieved by constructing a natural monoidal chain map
\begin{equation}\label{e:adams map 2}
	\theta_{\fX} \colon \cobar \sSchainsA(\fX,x) \to \cSchains(\loops_x \fX),
\end{equation}
for any pointed topological space $(\fX,x)$, and showing it to be a quasi-isomorphism for simply-connected spaces in \cite{adams1956cobar}, a hypothesis removed in \cite{rivera2018cubical} and, using a different argument, in \cite{rivera2019path}.
Here $\sSchainsA(\fX,x)$ denotes the coalgebra $\schainsA(\sSing(\fX,x))$ (using the notation introduced in \ref{ss:simplicial}) of (pointed) normalized singular chains.

We now describe the construction of Adams' map, whose combinatorial essence is the observation that the set of ascending chains in $0 < \dots < n$ containing $0$ and $n$, with the inclusion order, is isomorphic to the $(n-1)^\th$ cubical lattice.
To define Adams' map one uses a collection of continuous maps
\[
\set[\big]{\theta_n \colon \gcube^{n-1} \to P(\gsimplex^n;v_0,v_n)}_{n\in\N},
\]
satisfying for each $j$ the following conditions:
\begin{enumerate}
	\item $\theta_1(0)\colon [0,\sqrt{2}] \to \gsimplex^1$ is the path $\theta_1(0)(s) = v_0 + \frac{s}{\sqrt{2}}(v_1-v_0)$,
	\item $\theta_n \circ \delta_0^j = P(\delta^j) \circ \theta_{n-1}$, and
	\item $\theta_n \circ \delta_1^j =
	\big(P(\delta^{f^n_j}) \circ \theta_j\big) \cdot \big(P(\delta^{\ell^n_{n-j}}) \circ \theta_{n-j}\big)$,
\end{enumerate}
where, for $\epsilon \in \{0,1\}$,
\[
\delta_\epsilon^j = \id_{\gcube^{j-1}} \times \delta_\epsilon \times \id_{\gcube^{n-j-2}} \,\colon\ \gcube^{n-2} \to \gcube^{n-1},
\]
and
\[
\begin{split}
	\delta^{f^n_j} = \delta^n \circ \dotsb \circ \delta^{n-j+1} &\,\colon\ \gsimplex^j \to \gsimplex^n, \\
	\delta^{\ell^n_{n-j}} = \delta^{n-j-1} \circ \dotsb \circ \delta^{0} &\, \colon\ \gsimplex^{n-j} \to \gsimplex^n,
\end{split}
\]
are respectively the inclusions into the first and last faces of $\gsimplex^n$.

Adams showed the existence of a (non-unique) collection of such maps by induction, using the contractibility of $P(\gsimplex^n;v_0,v_n)$.
He then defined $\theta_{\fX}$ as follows.
For any singular $1$-simplex $\sigma \in \sSchainsA(\fX,x)$ let
\[
\theta_{\fX}[\sigma] = P(\sigma) \circ \theta_1 - c_x,
\]
where $(c_x \colon \gcube^0 \to \loops_x \fX) \in \cSchains(\loops_x \fX)$ is the singular $0$-cube determined by the constant loop (with $r=0$) at $x \in \fX$.
For any singular $n$-simplex $\sigma \in \sSchains(\fX,x)$ with $n>1$, let
\[
\theta_{\fX}[\sigma] = P(\sigma) \circ \theta_n.
\]
Since the underlying graded monoid structure of $\cobar \sSchainsA(\fX,x)$ is free, we may extend the above to a monoidal chain map $\theta_{\fX} \colon \cobar \sSchainsA(\fX,x) \to \cSchains(\loops_x \fX)$.
Conditions (1), (2), and (3) then imply that $\theta_{\fX}$ is a chain map.

\begin{remark}
	Adams originally worked with the sub-coalgebra $\sSchainsA(\fX,x)^1$ of $\sSchainsA(\fX,x)$ generated by singular simplices $\sigma \colon \gsimplex^n \to \fX$ collapsing the $1$-skeleton of $\gsimplex^n$ to $x \in \fX$.
	If $\fX$ is simply connected, the inclusion map $\sSchainsA(\fX,x)^1 \hookrightarrow \sSchainsA(\fX,x)$ induces a quasi-isomorphism on homology.
	In this case, the degree $1$ module in $\sSchainsA(\fX,x)^1$ is trivial, so that $\cobar \sSchainsA(\fX,x)^1$ is connected, i.e.
	isomorphic to the underlying ring $\k$ in degree~$0$.
	The cobar construction does not preserve quasi-isomorphisms in general.
	It does preserve quasi-isomorphisms between coaugmented coalgebras that are trivial in degree $1$.
\end{remark}

\subsection{An explicit choice}\label{explicitchoice}

We now construct an explicit collection of maps
\[
\set[\big]{\theta_n \colon \gcube^{n-1} \to P(\gsimplex^n;v_0,v_n)}_{n\in\N}
\]
satisfying the conditions in \cref{ss:adams maps}.

\begin{figure}
	\centering
	\includegraphics[scale=1]{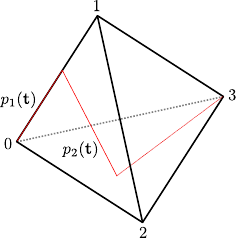}
	\caption{In red, the path $\theta_3(\bft) \in P(\gsimplex^3, v_0, v_3)$ associated to a $\bft \in \gcube^2$.}
	\label{f:theta1}
\end{figure}

Given $v,w \in \mathbb{R}^{n+1}$ denote by
\[
\gamma(v,w) \colon \big[0, \bars{w-v}\big] \to \mathbb{R}^{n+1}
\]
the straight line path from $v$ to $w$ parameterized by arc length, i.e.
\[
\gamma(v,w)(s) = v + \frac{s}{|w-v|}(w-v).
\]
For any $\bft=(t_1, \dots, t_{n-1}) \in \gcube^{n-1}$ we define $p_1(\bft), \dots, p_{n-1}(\bft)$ in $\gsimplex^n$ inductively by
\[
\begin{split}
	p_1(\bft) &= v_0+ t_1(v_1-v_0), \\
	p_j(\bft) &= p_{j-1}(\bft) + t_j(v_j-p_{j-1}(\bft)).
\end{split}
\]
We may now define
\[
\theta_n(\bft) \colon [0, r_{\bft}] \to \gsimplex^n,
\]
where
\[
r_{\bft} = \bars{p_1(\bft)-v_0} + \bars{p_2(\bft)-p_1(\bft)} + \dots + \bars{v_n-p_{n-1}(\bft)},
\]
as the piecewise linear path given by concatenating the straight line segments connecting the ordered sequence of points $v_0, p_1(\bft), \dots, p_n(\bft), v_n$, i.e.
\[
\theta_n(\bft) =
\gamma(v_0,p_1(\bft)) \cdot \gamma(p_1(\bft), p_2(\bft)) \cdot\ \dotsb\ \cdot \gamma(p_{n-1}(\bft),v_n).
\]
Please consult \cref{f:theta1} for an example illustrating this construction.

\begin{figure}[b]
	\centering
	\includegraphics[scale=.5]{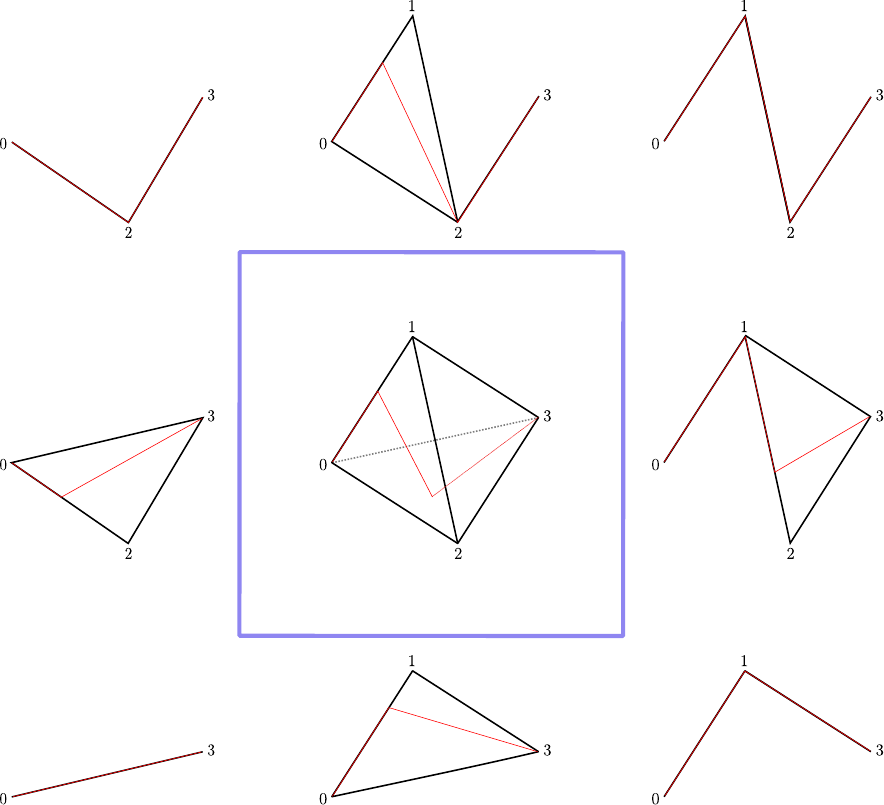}
	\caption{The faces of the $2$-cube in blue, labeled by a red element in their corresponding family of paths.}
	\label{f:theta3}
\end{figure}

A straightforward computation proves that the conditions in \cref{ss:adams maps} are satisfied by this collection.
A diagrammatical verification in low dimensions is provided in \cref{f:theta3}.

\begin{figure}
	\centering
	\includegraphics[scale=1]{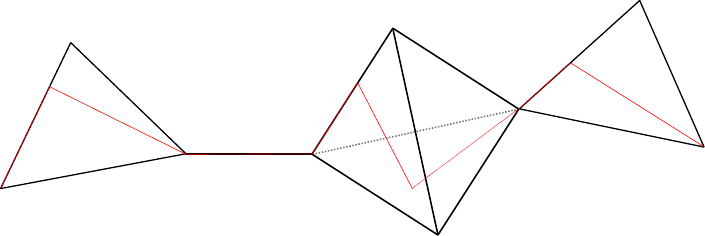}
	\caption{In red, the path $\theta_{(2,1,3,2)}(\bft)$ in $P(\gsimplex^2 \vee \gsimplex^1 \vee \gsimplex^3 \vee \gsimplex^2)$ associated to an element $\bft$ in $\gcube^{4}$.}
	\label{f:theta2}
\end{figure}

We can extend this collection $\{\theta_n \colon \gcube^{n-1} \to P(\gsimplex^n;v_0,v_n)\}$ to one of the form
\[
\set[\big]{\theta_{(n_1, \dots, n_k)} \colon \gcube^{n_1 + \dots + n_k-k} \to P(\gsimplex^{n_1} \vee \cdots \vee \gsimplex^{n_k})}
\]
where the \textit{topological necklace} $\gsimplex^{n_1} \vee \cdots \vee \gsimplex^{n_k}$ is obtained by identifying the last vertex of $\gsimplex^{n_i}$ with the first vertex of $\gsimplex^{n_{i+1}}$ for each $i = 1, \dots, k-1$ assuming each $n_i > 0$.
We do so by setting
\begin{multline*}
	\theta_{(n_1, \ldots, n_k)}(t_1, \ldots, t_{n_1+ \cdots + n_k-k}) \\ =
	\theta_{n_1}(t_1, \ldots, t_{n_1-1}) \cdot\ \dots\ \cdot \theta_{n_k}(t_{n_1+ \cdots n_{k-1}-k},\ldots,t_{n_1 + \cdots +n_k-k}).
\end{multline*}
Thus we can think of the topological necklace $\gsimplex^{n_1} \vee \cdots \vee \gsimplex^{n_k}$ as a space parameterizing a $(n_1+ \cdots + n_k-k)$-dimensional family of paths between the first and last vertices.
Please consult \cref{f:theta2} for an example illustrating this construction.

\subsection{Necklaces}

We now present a categorical viewpoint on the necklaces encountered in the previous subsection and their intimate relationship with Adams' cobar construction.
For any pointed space $(\fX,x)$, the underlying graded $\k$-module of $\cobar \sSchainsA(\fX,x)$ can be described as the graded $\k$-module freely generated by finite ordered sequences
\[
T = (\sigma_1, \dots, \sigma_k)
\]
of simplices $\sigma_i \in \sSing(\fX,x)$, with degree $\bars{T} = \bars{\sigma_1}+\dots+\bars{\sigma_k} - k$, modulo the sub $\k$-module generated by those sequences with at least one degenerate simplex.
The differential of $\cobar \sSchainsA(\fX,x)$ on a generator $T$ of degree $n$ can be expressed as a signed signed sum of all generators $T'$ of degree $n-1$. The simplices in the ordered sequence corresponding to each of these $T'$ are all faces of simplices in the ordered sequence corresponding to $T$. Note that two types of generators $T'$ appear in this differential: those that have the same length as $T$ and those that have exactly one more simplex. The first type corresponds to terms in the differential of $\sSchainsA(\fX,x)$ and the second type to terms in the Alexander-Whitney coproduct.
Furthermore, the monoid structure of $\cobar \sSchainsA(\fX,x)$ is induced by simply concatenating these ordered sequences of simplices.
This perspective suggests that $\cobar \sSchainsA(\fX,x)$ may be obtained by applying a normalized chains functor to certain cellular monoid, naturally associated to $(\fX,x)$, with cells labeled by finite ordered sequences of simplices in $\sSing(\fX,x)$.
A geometric construction reflecting this idea was described in \cite{baues1980geometry}.
We instead take a categorical approach and make this discussion precise through the framework of \textit{necklaces} and \textit{necklical sets} as we now discuss.

Consider the subcategory $\simplex_{*,*}$ of the simplex category $\simplex$ with the same objects and morphisms given by functors $f \colon [n] \to [m]$ satisfying $f(0) = 0$ and $f(n) = m$.
It is a strict monoidal category when equipped with the monoidal structure $[n] \ot [m] = [n+m]$, thought of as identifying the elements $n \in [n]$ and $0 \in [m]$, and unit given by $[0]$.
Heuristically, we may think of $\simplex_{*,*}$ as a category of models for cells parameterizing families of paths with fixed endpoints inside a simplex.

The \textit{necklace category} $\Nec$ is obtained from $\simplex_{*,*}$ as follows.
Thinking of $\simplex_{*,*}$ as a monoid in $\Cat$, we first apply the bar construction to it and produce a simplicial object in $\Mon_{\Cat}$ which, after realization, defines the strict monoidal category $\Nec$.
We denote the monoidal structure by
\[
\vee \colon \Nec \times \Nec \to \Nec.
\]
We describe $(\Nec, \vee)$ in more explicit terms.
The objects of $\Nec$, called \textit{necklaces}, are freely generated by the objects of $\simplex_{*,*}$ through the monoidal structure $\vee$.
Namely, the set of objects of $\Nec$ is the set of monomials
\[
\big\{ [n_1] \vee \dots \vee[n_k] \mid n_i, k \in \N_{>0} \big\}
\]
together with $[0]$ serving as the monoidal unit.

The morphisms of $\Nec$ are generated through the monoidal structure by the following four types of morphisms for all $n \in\N_{>0}$
\begin{enumerate}
	\item $\partial^j \colon [n-1] \to [n]$ for $j = 1, \dots, n-1$,
	\item $\Delta_{[j], [n-j]} \colon [j] \vee [n-j] \to [n]$ for $j = 1, \dots, n-1$,
	\item $\xi^j \colon [n+1] \to [n]$ for $j = 0, \dots, n$ and $n>0$, and
	\item $\xi^0 \colon [1] \to [0]$.
\end{enumerate}
We may identify $\Nec$ with a full sub-category of the category of double pointed simplicial sets $\sSet_{*,*}$ as follows.
Consider the functor
\[
\mathcal{S} \colon \Nec \to \sSet_{*,*}
\]
induced by sending any necklace $T = [n_1] \vee \dots \vee[n_k]$ to the simplicial set
\[
\mathcal{S}(T) = \simplex^{n_1} \vee \dots \vee \simplex^{n_k},
\]
where the wedge symbol now means we identify the last vertex of $\simplex^{n_i}$ with the first vertex of $\simplex^{n_{i+1}}$ for $i = 1, \dots, k-1$ and the two base points are given by the first vertex of $\simplex^{n_1}$ and the last vertex of $\simplex^{n_k}$.
Then $\mathcal{S}$ is a fully faithful functor, so $\Nec$ may be identified with the full sub-category of $\sSet_{*,*}$ having as objects those double pointed simplicial sets of the form $\simplex^{n_1} \vee \dots \vee \simplex^{n_k}$.
The \textit{dimension} of a necklace $T = [n_1] \vee\dots\vee [n_k]$ is defined by $\dim(T) = n_1 + \dots + n_k-k$.
Note there is a canonical homeomorphism \[\gsimplex^{n_1} \vee \cdots \vee \gsimplex^{n_k} \cong | \simplex^{n_1} \vee \cdots \vee \simplex^{n_k}|.\]

The category of \textit{necklical sets} $\Fun(\Nec^{\op}, \Set)$ becomes a (non-symmetric) monoidal category with monoidal structure induced from $(\Nec, \vee)$.
We denote the monoidal category of necklical sets by $\nSet$.

\begin{remark*}
	The category of necklaces was introduced in \cite{dugger2011rigidification} to give an explicit description of the homotopy coherent nerve functor and its left adjoint.
\end{remark*}

\subsection{From necklaces to cubes}

In \cref{explicitchoice} we described an explicit way of decomposing a topological necklace $\gsimplex^{n_1} \vee \cdots \vee \gsimplex^{n_k}$ into a family of paths connecting the first and last vertices parameterized by an $(n_1 + \cdots + n_k-k)$-dimensional cube.
This construction satisfies conditions (1), (2), and (3) in \cref{ss:adams maps}.
In particular, conditions (2) and (3) may be interpreted as saying that each face in the codimension $1$ boundary of such a cube of paths is in one-to-one correspondence with codimension $1$ ``sub-necklaces" inside $\gsimplex^{n_1} \vee \cdots \vee \gsimplex^{n_k}$ connecting the first and last vertices.
Furthermore, sub-necklaces in $\gsimplex^{n_1} \vee \cdots \vee \gsimplex^{n_k}$ are in one-to-one correspondence with the poset
\[
\set[\Big]{J \subseteq \set{0,\dots,n_1+\cdots+n_k-k} \mid 0,n_1+\cdots+n_k-k \in J}
\]
ordered by inclusion, which has a canonical cubical structure.
We build upon this observation to describe a functorial relation between necklaces and cubes.
This will be used to explain how Adams' map is induced by a deeper categorical construction.

We begin by defining a monoidal functor
\[
\cP \colon \Nec \to \cube
\]
as follows.
First define $\cP[0] = 2^0$.
On any other necklace $T \in \Nec$, define $\cP(T) = 2^{\dim(T)}$.
In order to define $\cP$ on morphisms, it is sufficient to consider the following cases.
\begin{enumerate}
	\item For any coface map $\partial^j \colon [n] \to [n+1]$ such that $0< j<{n+1}$, define $\cP(\partial^j) \colon 2^{n-1}\to 2^{n}$ to be the cubical coface functor $\cP(f)= \delta_0^{j}$.

	\item For any $\Delta_{[j], [n+1-j]} \colon [j] \vee [n+1-j] \to [n+1]$ such that $0<j<n+1$, define
	\[
	\cP(\Delta_{[j], [n+1-j]}) \colon 2^{n-1}\to 2^{n}
	\]
	to be the cubical coface functor $\cP(f)=\delta_1^{j}$.

	\item We now consider codegeneracy maps of the form $\xi^j \colon [n+1] \to [n]$ for $n>0$.
	If $j=0$ or $j=n$, define $\cP(f) \colon 2^n \to 2^{n-1}$ to be the cubical codegeneracy functor $\cP(s^j)= \varepsilon^{j}$.
	If $0<j<n$, define $\cP(s^j) \colon 2^n \to 2^{n-1}$ to be the cubical coconnection functor $\gamma^{j}$.

	\item For $\xi^0 \colon [1] \to [0]$ define $\cP(\xi^0) \colon 2^0 \to 2^0$ to be the identity functor.
\end{enumerate}

\begin{remark*}
	The functor $\cP$ is neither faithful or full.
	However, for any necklace $T' \in \Nec$ with $\dim(T')=n+1$ and any cubical coface functor $\delta_{\epsilon}^j \colon 2^n \to 2^{n+1}$ for $0 \leq j \leq n+1$, there exists an map $f \colon T \hookrightarrow T'$, where $T \in \Nec$ with $\dim(T)=n$ such that $\mathcal{S}(f) \colon \mathcal{S}(T) \hookrightarrow \mathcal{S}(T')$ is an injective morphism in $\Nec$ and $\cP(f) = \delta_{\epsilon}^j$.
\end{remark*}

The functor $\cP \colon \Nec \to \cube$ induces an adjunction between $\cSet$ and $\nSet$ with right and left adjoint functors given respectively by
\[
\cP^\ast \colon \cSet \to \nSet,
\qquad \text{and} \qquad
\cP_{!} \colon \nSet \to \cSet.
\]
Explicitly, for a cubical set $Y \colon \cube^\op \to \Set$,
\[
\cP^\ast(Y) = Y \circ \cP^\op,
\]
and for a necklical set $K \colon \Nec^\op \to \Set$,
\[
\cP_{!}(K) \ =
\colim_{\yoneda(T) \to K} \cP(T) \ \cong
\colim_{\yoneda(T) \to K} \cube^{\dim(T)},
\]
where $\yoneda \colon \Nec \to \nSet$ is the Yoneda embedding.
Since $\cP$ is a monoidal functor, $\cP_{!} \colon \nSet \to \cSet$ is monoidal as well.

\subsection{Cubical cobar construction}\label{ss:cubical cobar}

Using the framework of necklical sets, we may reinterpret Baues' geometric cobar construction \cite{baues1980geometry} as a functor
\[
\ncobar \colon \sSet^0 \to \Mon_{\nSet},
\]
which we now define.

For any reduced simplicial set $X$, we define a necklical set $\ncobar(X) \colon \Nec^\op \to \Set$ having as necklical cells all necklaces inside $X$; namely
\[
\ncobar(X) \, = \! \colim_{\mathcal{S}(T) \to X} \yoneda(T).
\]
The monoidal structure $\vee \colon \Nec \times \Nec \to \Nec$ given by concatenation of necklaces induces a natural product
\[
\ncobar(X) \ot \ncobar(X) \to \ncobar(X)
\]
making $\ncobar(X)$ into a monoid in $\nSet$.

We may now define the \textit{cubical cobar construction}
\[
\ccobar \colon \sSet^0 \to \Mon_{\cSet}
\]
as the composition
\[
\ccobar = \cP_! \, \circ \ncobar.
\]
Since $\cP_!$ is monoidal, $\ccobar(X)$ is a monoid in $\cSet$.

\begin{remark}
	This reinterpretation of Baues' construction in terms of cubical sets was also studied in \cite{rivera2018cubical}.
	In this reference, it is also proven that the composition of functor $\mathcal{T} \circ \ccobar$, where $\mathcal{T} \colon \Mon_{\cSet} \to \Mon_{\sSet}$ is the triangulation functor, coincides with the left adjoint of the homotopy coherent nerve functor restricted to $\sSet^0.$
\end{remark}

\subsection{Relation to the cobar construction}

We now relate the cubical cobar functor $\ccobar \colon \sSet^0 \to \Mon_{\cSet}$ to the cobar construction $\cobar \colon \coAlg^\ast \to \Mon_{\Ch}$ (\cref{ss:cobar construction}).

\begin{theorem}\label{t:ccobar and cobar}
	There is a natural isomorphisms of functors
	\[
	\cchains \ccobar \cong \cobar \schainsA \colon \sSet^0 \to \Mon_{\Ch}.
	\]
\end{theorem}

\begin{proof}
	Denote by $\iota_n \in (\cube^n)_n$ the top dimensional non-degenerate element of the standard $n$-cube $\cube^n$.
	Note that for a reduced simplicial set $X$, we may represent any non-degenerate $n$-cube $\alpha \in (\cP_!(\ncobar(X)))_n$ as a pair $\alpha = [\sigma \colon \yoneda(T) \to X, \iota_n]$ for some $T = [n_1] \vee \dots \vee [n_k] \in \Nec$ with $\dim(T) = n_1 + \dots + n_k - k = n$.

	To define a monoidal chain map
	\[
	\varphi_X \colon \cchains(\cP_!(\ncobar(X))) \xra{\cong} \cobar \schainsA(X)
	\]
	it suffices to define it on any generator of the form $\alpha=[\sigma \colon \simplex^{n+1} \to X, \iota_{n}]$, i.e., when $T$ is of the form $T = [n+1]$, for some $n \geq 0$.
	If $n = 0$ let $\varphi_X(\alpha)= [\overline{\sigma}] + 1_\k$, where $[\overline{\sigma}] \in \s{-1} \overline{ \schains(X)} \subset \cobar \schainsA(X)$ denotes the (length $1$) generator in the cobar construction of $\schains(X)$ determined by $\sigma \in X_{n+1}$ and $1_\k$ denotes the unit of the underlying ring $\k$.
	If $n > 0$, we let $\varphi_X(\alpha)=[\overline{\sigma}]$.
	A straightforward computation yields that this gives rise to a well defined isomorphism of algebras, which is compatible with the differentials, and natural with respect to maps of simplicial sets.
\end{proof}

A similar result to \cref{t:ccobar and cobar} was observed in the case of $1$-reduced simplicial sets in \cite[Section~3.5]{berger1995loops}.

\subsection{Factorization of Adams's map}\label{ss:factorization of adams}

Adams's comparison map can be factored as a composition
\begin{equation}\label{e:factorization of adams}
	\theta_\fX \colon \cobar \sSchainsA(\fX,x) \xrightarrow{\cong}
	\cchains \ccobar(\sSing(\fX,x)) \xrightarrow{\cchains(\Theta)}
	\cSchains(\loops_x \fX).
\end{equation}
The first map is the monoidal isomorphism induced by \cref{t:ccobar and cobar}.
The second map is given by applying chains $\cchains \colon \Mon_{\cSet} \to \Mon_{\Ch}$ to the map of monoidal cubical sets
\[
\Theta \colon \ccobar(\sSing(\fX,x)) \to \cSing(\loops_x \fX)
\]
determined through the monoid structure by sending an $n$-simplex $(\sigma \colon \gsimplex^n \to \fX)$ to the singular $(n-1)$-cube
\[
P(\sigma) \circ \theta_n \colon \gcube^{n-1} \to \loops_x \fX,
\]
where the maps $\theta_n \colon \gcube^{n-1} \to P(\gsimplex^n;0,n)$ are discussed in \cref{ss:adams maps}.

\subsection{A monoidal coalgebra structure on the cobar construction}

We follow \cite{baues1998hopf} to construct a coalgebra structure on $\cobar \sSchainsA(\fX,x)$, compatible with both the differential and monoid structure, such that Adams's map becomes a map of monoids in the category of coalgebras.

Recall the Serre coalgebra lift $\cchainsA \colon \cSet \to \coAlg$ of $\cchains \colon \cSet \to \Ch$, the unique monoidal functor defined by the coalgebra structure on $\chains(\cube^1)$.
Since $\cchainsA$ is monoidal and $\ccobar(\sSing(\fX,x))$ is a monoid, $\cchainsA(\ccobar(\sSing(\fX,x)))$ is a monoid in $\coAlg$.
Similarly, the lift $\cchainsA$ equips $\cSing(\loops_x \fX)$ with a natural monoidal coalgebra structure as well.
Consequently, the isomorphism in (\ref{e:factorization of adams}) endows $\cobar \sSchainsA(\fX,x)$ with a natural monoidal coalgebra structure making $\theta_\fX$ into natural map of monoidal coalgebras.

\section{Monoidal \pdfEinfty-structures}

In this section we recall the model $\UM$ of the $E_\infty$-operad, whose category of coalgebras we show to be monoidal.
We use this structure to construct a monoidal functor $\cchainsUM \colon \cSet \to \coAlg_\UM$ extending the Serre coalgebra structure.
This endows for any pointed topological space both $\cobar \sSchainsA(\fX,x)$ and $\cSchainsA(\loops_x \fX)$ with the structure of a monoidal $E_\infty$-coalgebra which is preserved by Adams' map.

\subsection{\pdfEinfty-operads}

Recall that operads control algebraic structures with either one input and multiple outputs or vice-versa.
In this article, we work with dg operads, i.e. operads in the monoidal category $(\Ch, \otimes, \mathbb{1})$; we refer to \cite{vallette2012operads} for more details.

For example, coalgebras, as defined in \cref{ss:coalgebras}, are controlled by the operad $\As$ generated by two elements in degree $0$
\[
\coproduct\,, \quad \counit\,,
\]
modulo the relations
\[
\leftcounitality\,, \quad \rightcounitality\,, \quad \coassociativity\,.
\]
Let $\sigma \in \sym_2$ be the non-identity transposition.
The operad $\Com$ controlling cocommutative coalgebras is obtained by adding the relation
\[
\cocommutativity
\]
to this presentation.
We are interested in compatibly resolving the (trivial) symmetric group actions on $\Com$ associated to the permutation of factors.
An $E_\infty$-operad is an operad $\cO$ quasi-isomorphic to $\Com$ for which the action of $\sym_r$ on $\cO(r)$ is free for each $r \in \N$.
An $E_\infty$-coalgebra structure on a chain complex $C$ is an operad morphism
\[
\cO \to \coEnd(C) = \set{\Hom(C, C^{\ot r})}_{r \in \N}.
\]

\subsection{A finitely presented \pdfEinfty-prop}

A prop is an object controlling algebraic structures with multiple inputs and outputs.
We refer to \cite{markl2008props} for a detailed exposition.
We recall the following construction from \cite{medina2020prop1}.
The prop $\cM$ is generated by adding to the presentation of $\As$ a generator and a relation.
More specifically, a generator in degree $1$ with boundary
\[
\product \ \raisebox{3pt}{$\xmapsto{\bd}$} \ \boundary\,,
\]
and the relation
\[
\qquad \productcounit\,.
\]
The importance of this construction is that the operad $\UM = \set{\cM(1,r)}_{r\in\N}$ obtained by restriction of structure is an $E_\infty$-operad.

\subsection{\pdfM-bialgebras}

An $\cM$-bialgebra structure on $C$ is a prop morphism
\[
\cM \to \biEnd(C) = \set{\Hom(C^{\ot s}, C^{\ot r})}_{r,s\in\N}.
\]
More explicitly, is a coalgebra $(C, \Delta, \aug)$ together with a degree~1 product satisfying for any $a,b \in C$ that:
\begin{align}
	\label{eq:M-bialg def 1}
	\aug(a \ast b) =\ & 0, \\
	\label{eq:M-bialg def 2}
	\bd (a \ast b) =\ & \bd a \ast b - (-1)^{a} a \ast \bd b + \aug(a) b - (-1)^{a} a \aug(b).
\end{align}
The signs appearing in Identity \eqref{eq:M-bialg def 2} are a result of the Koszul sign convention.

Any $\cM$-bialgebra structure induces an $E_\infty$-coalgebra structure.
More explicitly, let $(C, \Delta, \aug, \ast)$ be an $\cM$-bialgebra.
The collection of all maps $\set{C \to C^{\ot r}}_{r \in \N}$ generated by $\Delta$, $\aug$ and $\ast$ makes $C$ into an $E_\infty$-coalgebra, specifically into an $\UM$-coalgebra.

\subsection{\pdfEinfty-structure on simplicial chains}

We recall the construction of an $E_\infty$-extension of the Alexander--Whitney coalgebra structure on simplicial chains introduced in \cite{medina2020prop1}.
Let us start by considering the representable simplicial sets.
The coalgebra $\chains(\simplex^n)$ can be made into a natural $\cM$-bialgebra considering an algebraic version of the \textit{join product} defined by
\begin{equation*}
	\left[v_0, \dots, v_p \right] \ast \left[v_{p+1}, \dots, v_q\right] =
	\begin{cases}
		(-1)^{p} \sign(\pi) \left[v_{\pi(0)}, \dots, v_{\pi(q)}\right] &
		\text{ if } v_i \neq v_j \text{ for } i \neq j, \\
		\hfil 0 & \text{ if not},
	\end{cases}
\end{equation*}
where $\pi$ is the permutation that orders the vertices.
A Kan extension of the induced $\UM$-coalgebra structure on $\chains(\simplex^n)$ defines a lift $\schainsUM$ of the Alexander--Whitney coalgebra to the category of $E_\infty$-coalgebras defined by $\UM$.
That is to say, these functors fit in a commutative diagram
\begin{equation*}
	\begin{tikzcd}[column sep=normal, row sep=small]
		& \coAlg_{\UM} \arrow[d] \\
		\sSet \arrow[r]
		\arrow[ur,out=60, in=180, "\schainsUM"]
		\arrow[r, "\schainsA"]
		& \coAlg,
	\end{tikzcd}
\end{equation*}
where the vertical arrow is the obvious forgetful functor.

We remark that this $E_\infty$-structure generalizes those previously introduced in \cite{mcclure2003multivariable,berger2004combinatorial} in the sense that any cooperation $\chains \to \chains^{\ot r}$ arising from the action of the Barratt-Eccles or surjection operad can be expressed as a cooperation arising from the action of $\UM$.

\subsection{\pdfEinfty-structure on cubical chains}\label{ss:cubical e-infty}

We recall the construction of an $E_\infty$-extension of the Serre coalgebra structure on cubical chains introduced in \cite{medina2022cube_einfty}.
Let us first consider representable cubical sets.
The coalgebra structure on $\chains(\cube^n)$ can be made into a natural $\cM$-bialgebra considering a product defined using the following notation.
For a basis element $x = x_1 \ot \dotsb \ot x_n$ of $\chains(\cube^n)$ and an integer $\ell \in \set{1,\dots,n}$ we write
\begin{align*}
	x_{<\ell} & = x_1 \ot \dotsb \ot x_{\ell-1}, \\
	x_{>\ell} & = x_{\ell+1} \ot \dotsb \ot x_n,
\end{align*}
with the convention $x_{<1} = x_{>n} = 1 \in \Z$.
Then, for two such basis elements $x$ and $y$ we set
\begin{equation}\label{eq:product on cubes}
	(x_1 \ot \dotsb \ot x_n) \ast (y_1 \ot \dotsb \ot y_n) =
	\sum_{i=1}^n x_{<i}\, \epsilon(y_{<i}) \ot (x_i \ast y_i) \ot \epsilon(x_{>i}) \, y_{>i},
\end{equation}
where the only non-zero values of $x_i \ast y_i$ are
\[
[0] \ast [1] = [0, 1], \qquad [1] \ast [0] = -[0, 1].
\]
A Kan extension of the induced $\UM$-coalgebra structure on $\chains(\cube^n)$ defines a lift $\cchainsUM$ of the Serre coalgebra structure to the category of $E_\infty$-coalgebras defined by $\UM$.
That is, a commutative diagram
\begin{equation}\label{eq:lift to e-infty cubical}
	\begin{tikzcd}[column sep=normal, row sep=small]
		& \coAlg_{\UM} \arrow[d] \\
		\cSet \arrow[r]
		\arrow[ur,out=60, in=180, "\cchainsUM"]
		\arrow[r, "\cchainsA"]
		& \coAlg.
	\end{tikzcd}
\end{equation}

\begin{remark*}
	The product defined above in \cref{eq:product on cubes} differs from the one defined in \cite{medina2022cube_einfty} by the sign $(-1)^x$.
	The convention used here is more natural as we will see in \cref{ss:cube_einfty revisited}.
\end{remark*}

\subsection{Monoidal structure}\label{ss:monoidal structures}

In this subsection we describe an extension of the tensor product of coalgebras to $\cM$-bialgebras and $\UM$-coalgebras.

\begin{lemma}\label{l:monoidal M-bialg}
	Let $C$ and $C'$ be $\M$-bialgebras.
	The coalgebra $C \ot C'$ is a natural $\M$-bialgebra with
	\begin{equation}\label{eq:monoidal_product}
		(a \ot b) \ast (c \ot d) =
		a \aug(c) \ot (b \ast d) + (a \ast c) \ot \aug(b) d,
	\end{equation}
	for any $a,c \in C$ and $b,d \in C'$.
\end{lemma}

\begin{proof}
	We verify Identity \eqref{eq:M-bialg def 1} using that $\aug(b \ast d) = \aug(a \ast c) = 0$,
	\begin{align*}
		\aug\big((a \ot b) \ast (c \ot d)\big) =&
		\aug\big(a \aug(c) \ot (b \ast d)\big) + \aug\big((a \ast c) \ot \aug(b) d\big) \\ =&
		\aug(a) \aug(c) \ot \aug(b \ast d) + \aug(a \ast c) \ot \aug(b) \aug(d) \\ =& \ 0.
	\end{align*}
	To verify Identity \eqref{eq:M-bialg def 2} we need to show that
	\begin{equation}\label{eq:monoidal appendix}
		\begin{split}
			\bd \big((a \ot b) \ast (c \ot d)\big) =\, &
			\bd (a \ot b) \ast (c \ot d) - (-1)^{a+b} (a \ot b) \ast \bd (c \ot d) \\ +\, &
			\aug(a \ot b) (c \ot b) - (-1)^{a+b} (a \ot b) \aug(c \ot d).
		\end{split}
	\end{equation}
	Let us start computing the left hand side of the above expression.
	\begin{align*}
		\bd\big((a \ot b) \ast (c \ot d)\big) =\ &
		\bd\big(a \aug(c) \ot (b \ast d) + (a \ast c) \ot \aug(b) d\big) \\ =\ &
		\bd a \aug(c) \ot (b \ast d) \,+\, (-1)^{a} a \aug(c) \ot \bd(b \ast d) \\ +\ &
		\bd (a \ast c) \ot \aug(b) d \,+\, (-1)^{a+c+1} (a \ast c) \ot \aug(b) \bd d.
	\end{align*}
	Using that $C$ and $C'$ satisfy Identity \eqref{eq:M-bialg def 2} we have that
	\begin{align*}
		\bd\big((a \ot b) \ast (c \ot d)\big) =\ &
		\bd a \aug(c) \ot (b \ast d) \\ +\ &
		(-1)^a a \aug(c) \ot \big(\bd b \ast d + (-1)^{b+1} b \ast \bd d + \aug(b) d + (-1)^{b+1} b \aug(d)\big) \\ +\ &
		\big(\bd a \ast c + (-1)^{a+1} a \ast \bd c + \aug(a) c + (-1)^{a+1} a \aug(c)\big) \ot \aug(b) d \\ +\ &
		(-1)^{a+c+1} (a \ast c) \ot \aug(b) \bd d.
	\end{align*}
	Inspecting this expression we label terms which sum to $\bd\big((a \ot b) \ast (c \ot d)\big)$:\vspace*{-10pt}
	\begin{minipage}[t]{0.5\textwidth}
		\begin{align}& \label{x4}
			(-1)^0 \bd a \aug(c) \ot (b \ast d) \\& \label{x5}
			(-1)^{a} a \aug(c) \ot (\bd b \ast d) \\& \label{x6}
			(-1)^{a+b+1} a \aug(c) \ot (b \ast \bd d) \\& \label{x7}
			(-1)^a a \aug(c) \ot \aug(b) d \\& \label{x8}
			(-1)^{a+b+1} a \aug(c) \ot b \aug(d)
		\end{align}
	\end{minipage}
	\begin{minipage}[t]{0.5\textwidth}
		\begin{align}& \label{x9}
			(-1)^0 (\bd a \ast c) \ot \aug(b) d \\& \label{x10}
			(-1)^{a+1} (a \ast \bd c) \ot \aug(b) d \\& \label{x11}
			(-1)^0 \aug(a) c \ot \aug(b) d \\& \label{x12}
			(-1)^{a+1} a \aug(c) \ot \aug(b) d \\& \label{x13}
			(-1)^{a+c+1} (a \ast c) \ot \aug(b) \bd d \quad ||
		\end{align}
		\vskip\parskip
	\end{minipage}
	Additionally, since
	\begin{equation*}
		\bd (a \ot b) \ast (c \ot d) \,=\,
		(\bd a \ot b) \ast (c \ot d) \,+\,
		(-1)^a (a \ot \bd b) \ast (c \ot d),
	\end{equation*}
	the following labeled terms sum to $\bd (a \ot b) \ast (c \ot d)$:
	\begin{align}& \label{x14}
		(-1)^0 \bd a \aug(c) \ot (b \ast d) \\& \label{x15}
		(-1)^0 (\bd a \ast c) \ot \aug(b) d \\& \label{x16}
		(-1)^{a} a \aug(c) \ot (\bd b \ast d) \quad ||
	\end{align}
	Similarly, since
	\begin{equation*}
		(a \ot b) \ast \bd (c \ot d) \,=\,
		(a \ot b) \ast (\bd c \ot d) \,+\,
		(-1)^{c} (a \ot b) \ast (c \ot \bd d),
	\end{equation*}
	the following labeled terms sum to $(-1)^{a+b+1}(a \ot b) \ast \bd (c \ot d)$:
	\begin{align}& \label{x17}
		(-1)^{a+1} (a \ast \bd c) \ot \aug(b) d \\& \label{x18}
		(-1)^{a+b+1} a \aug(c) \ot (b \ast \bd d) \\& \label{x19}
		(-1)^{a+c+1} (a \ast c) \ot \aug(b) \bd d \quad ||
	\end{align}
	We have the following matching pairs of labeled summands.
	\begin{center}
		\eqref{x4}-\eqref{x14} :
		\eqref{x5}-\eqref{x16} :
		\eqref{x6}-\eqref{x18} :
		\eqref{x7}-\eqref{x12} :
		\eqref{x9}-\eqref{x15} :
		\eqref{x10}-\eqref{x17} :
		\eqref{x13}-\eqref{x19}.
	\end{center}
	Additionally, the sum of the unmatched terms \eqref{x11} and \eqref{x8} correspond to
	\[
	\aug(a \ot b) (c \ot b) - (-1)^{a+b} (a \ot b) \aug(c \ot d),
	\]
	which concludes the verification of both Identity \eqref{eq:monoidal appendix} and the lemma.
\end{proof}

We now give a more conceptual description of the monoidal structure on $\biAlg_\M$, which generalizes to $\coAlg_\UM$.
The prop $\cM$ is obtained by applying the functor of cellular chains to a CW prop \cite{medina2021prop2}.
These cells are in fact cubical, generated through compositions by the generators
\[
\stackrel{\cellcoproduct}{\ \gcube^0,} \qquad\qquad
\stackrel{\cellcounit}{\ \gcube^0,} \qquad\qquad
\stackrel{\cellcoboundary}{\ \gcube^1.}
\]
The Serre diagonal defines a diagonal $\Delta_\cM$ on $\cM$ compatible with its prop structure.
More specifically, for any basis element $\mu$ in $\M$ we have a chain map $\phi_\mu \colon \chains(\cube^n) \to \cM$ with $\phi_\mu(2^n) = \mu$, and $\Delta_\cM(\mu)$ is set to be $\phi_D^{\ot 2} \circ \Delta (2^n)$.
Crucially, $\Delta_\cM$ acts on the generators by
\begin{align*}
	\coproduct \quad &\mapsto \quad \coproduct \ \raisebox{3pt}{$\ot$} \; \coproduct \,, \\
	\counit \ \quad &\mapsto \ \quad \counit \ \raisebox{3pt}{$\ot$} \; \counit \,, \\
	\product \quad &\mapsto \, \quad \rightboundary \ \raisebox{3pt}{$\ot$} \product \ \raisebox{3pt}{$+$} \ \product \raisebox{3pt}{$\ot$} \ \leftboundary \,,
\end{align*}
which recovers the statement of \cref{l:monoidal M-bialg}.
The structure preserving diagonal on the prop $\cM$ induces one on the operad $\UM$ and defines, as usual for so-called Hopf operads, a monoidal structure on $\coAlg_\UM$.

\subsection{Revisiting the \pdfEinfty-coalgebra structure on cubical chains}\label{ss:cube_einfty revisited}

We show that the monoidal structure on $\UM$-coalgebras recover the $E_\infty$-structure on cubical chains defined in \cref{ss:cubical e-infty}.
In fact, we have the following stronger statement at the level of $\M$-bialgebras.

\begin{theorem}\label{t:cubical e-infty chains are monoidal}
	For any $n \in \N$, the $\M$-bialgebra structure on $\chains(\cube^n)$ agrees with the monoidal extension of the $\M$-bialgebra structure on $\chains(\cube^1)$.
\end{theorem}

\begin{proof}
	Since the coalgebra part agrees by definition, we focus on the product.
	We will proceed by induction with the base case holding trivially.
	Let $x = x_1 \ot \dotsb \ot x_n$ and $y = y_1 \ot \dotsb \ot y_n$ be two elements in $\chains(\cube^n)$.
	The following straightforward computation verifies that $x \ast y$ in $\chains(\cube^n)$ corresponds to $(x_{<n} \ot x_n) \ast (y_{<n} \ot y_n)$ in $\chains(\cube^{n-1}) \ot \chains(\cube^1)$:
	\begin{align*}
		x \ast y &=
		\sum_{i=1}^n x_{<i}\, \epsilon(y_{<i}) \ot x_i \ast y_i \ot \epsilon(x_{>i}) \, y_{>i} \\ &=
		\sum_{i=1}^{n-1} x_{<i}\, \epsilon(y_{<i}) \ot x_i \ast y_i \ot \epsilon(x_{>i}) \, y_{>i} \ +\
		x_{<n}\, \epsilon(y_{<n}) \ot x_n \ast y_n \\ &=
		x_{<n} \ast y_{<n} \ot \aug(x_n)\, y_n \ +\ x_{<n}\, \epsilon(y_{<n}) \ot x_n \ast y_n \\ &=
		(x_{<n} \ot x_n) \ast (y_{<n} \ot y_n),
	\end{align*}
	which concludes the proof.
\end{proof}

\subsection{A monoidal \pdfEinfty-coalgebra structure on the cobar construction}\label{ss:e-infty on cobar}

Using \cref{t:cubical e-infty chains are monoidal} and the natural equivalence of functors $\cchains \ccobar \cong \cobar \schainsA$ proven in \cref{t:ccobar and cobar}, we can transfer a monoidal $E_\infty$-structure to the Adams' cobar construction of any reduced simplicial set, extending the monoidal coalgebra structure defined by Baues.
This is summarized by the following diagram commuting up to isomorphisms:
\[
\begin{tikzcd}
	& \Mon_{\coAlg_\UM} \arrow[d] \\
	\Mon_{\cSet} \arrow[ru, "\cchainsUM", out=70, in=180] \arrow[r, "\cchainsA"]
	& \Mon_{\coAlg} \arrow[d] \\
	\sSet^0 \arrow[r, "\cobar \schainsA"] \arrow[u, "\ccobar"]
	& \Mon_{\Ch}.
\end{tikzcd}
\]
Furthermore, using the factorization of Adams' map described in \cref{ss:factorization of adams} and the naturality of our monoidal $E_\infty$-structure, we conclude that for any pointed topological space $\fX$,
\[
\theta_\fX \colon \cobar \sSchainsA(\fX,x) \to \cSchainsUM(\loops_x \fX)
\]
is a monoidal quasi-isomorphism of $E_{\infty}$-coalgebras, as announced in the introduction.
	\sloppy
	\printbibliography
\end{document}